\documentclass[a4paper]{jpconf}
\usepackage{graphicx}
\usepackage{latexsym}
\usepackage{graphics}
\usepackage{amssymb}
\usepackage{amsmath}
\usepackage{amsthm}
\usepackage{mathrsfs}
\usepackage{media9}
\usepackage{multimedia}
\usepackage{fancybox}
\usepackage{minibox}
\usepackage[usestackEOL]{stackengine}
\usepackage{fancyhdr}
\usepackage{booktabs}
\usepackage{lipsum}
\usepackage{pgf,tikz}
\usetikzlibrary{arrows}
\usetikzlibrary[patterns]
\usepackage{capt-of}
\usepackage{caption}
\usepackage{float}
\usepackage{color}
\usepackage{float}
\usepackage{animate}
\usepackage{graphicx}

\definecolor{qqqqcc}{rgb}{0,0,0.8}
\definecolor{ffqqqq}{rgb}{1,0,0}
\definecolor{rvwvcq}{rgb}{0.08235294117647059,0.396078431372549,0.7529411764705882}

\usepackage{braids}

\usepackage{mathtools}
\mathtoolsset{showonlyrefs=true}

\newtheorem*{bellingeri}{Proposition ([2, Theorem 1.2])}
\newtheorem*{GR}{Grauert-Remmert extension theorem [3, XII.5.4]}
\newtheorem*{GAGA}{GAGA principle [4]}
\newtheorem*{proj-ext}{Corollary (Analytic to algebraic extension of coverings)}
\newtheorem*{monodromy}{Corollary (Branched coverings via monodromy)}
\newtheorem*{theoremP1}{Theorem [1, Theorem 1]}
\newtheorem*{theoremP2}{Theorem [1, Theorem 2]}

\theoremstyle{definition}
\newtheorem*{definition-conf}{Definition (Configuration spaces)}
\newtheorem*{definition-braid}{Definition (Surface braid groups)}
\newtheorem*{definition-pluri}{Definition (Pluricanonical maps)}
\newtheorem*{definition-alb}{Definition (Albanese map)}
\newtheorem*{definition-gen-cov}{Definition (Generic covering)}
\newtheorem*{definition-gen-mod}{Definition (Generic monodromy representation)}

\newtheorem*{remark}{Remark}
\newtheorem*{example-3}{The case $n=3$}

\newtheorem*{example-2}{The case $n=2$}

\renewcommand{\SS}{\mathrm{Sym^2}(C_2)}
\renewcommand{\to}{\longrightarrow}

\begin{document}
\title{Representations of braid groups and construction of projective surfaces}

\author{Francesco Polizzi}

\address{Dipartimento di Matematica e Informatica \\
Universit\`a della Calabria (Italy)}

\ead{polizzi@mat.unical.it}

\begin{abstract}
Braid groups are an important and flexible tool used in several areas of science, such as Knot Theory (Alexander's theorem), Mathematical Physics (Yang-Baxter's equation) and Algebraic Geometry (monodromy invariants). In this note we will focus on their algebraic-geometric aspects, explaining how the representation theory of higher genus braid groups can be used to produce interesting examples of projective surfaces defined over the field of complex numbers.
\end{abstract}

\section{Introduction}
The classification of compact, complex surfaces $S$ of general type with $\chi(\mathcal{O}_S)=1$, i.e. $p_g(S)=q(S)$, is currently an active area of research. For all these surfaces we have $p_g \leq 4$, and the cases $p_g=q=4$ and $p_g=q=3$ are nowadays completely described. Regarding the case $p_g=q=2$, a complete classification has been recently obtained when $K_S^2=4$: in fact, these are surfaces on the Severi line $K_S^2=4 \chi (\mathcal{O}_S)$. By contrast, the classification in the case $p_g=q=2$, $K_S^2 \geq 5$ is still missing, albeit some interesting examples were recently discovered. We refer the reader to the paper \cite{1} and the references contained therein for a historical account on the subject and more details. 

The purpose of this note is to show how monodromy representations of braid groups can be concretely applied to the fine classification of surfaces with $p_g=q=2$ and maximal Albanese dimension, allowing one to rediscover old examples and to find new ones.

The idea is to consider degree $n$, generic covers of $\SS$, the symmetric square of a smooth curve of genus $2$, simply branched over the diagonal $\delta$. In fact, if such a cover exists, then it is a smooth surface $S$ with
\begin{equation*}
\chi(\mathcal{O}_S)=1, \quad K_S^2=10-n.
\end{equation*} 

On the other hand, by the Grauert-Remmert extension theorem and the GAGA principle, isomorphism classes of degree $n$, connected covers 
\begin{equation*}
f \colon S \to \SS, 
\end{equation*}
branched at most over $\delta$, correspond to group homomorphisms 
\begin{equation*}
\varphi \colon \pi_1(\SS - \delta) \to \mathfrak{S}_n
\end{equation*}
with transitive image, up to conjugacy in $\mathfrak{S}_n$. The group $\pi_1(\SS - \delta)$ is isomorphic to $\mathsf{B}_2(C_2)$, the braid group on two strings on $C_2$; furthermore, our condition that the branching is simple can be translated by requiring that $\varphi(\sigma)$ is a transposition, where $\sigma$ denotes the homotopy class in $\SS - \delta$ of a topological loop in $\SS$ that ``winds once around $\delta$". 

A group homomorphism $\mathsf{B}_2(C_2) \to \mathfrak{S}_n$ satisfying the requirements above will be called a \emph{generic monodromy representation} of $\mathsf{B}_2(C_2)$. By using the Computer Algebra System \verb|GAP4|  we computed the number of generic monodromy representations  for $2 \leq n \leq 9$. In particular, such a number is zero for $n \in \{5, \, 7, \, 9 \}$, so there exist no generic covers in these cases.

As an application of the general theory, we end the paper with a detailed discussion of the cases $n=2$ and $n=3$.

\section{Braid groups on closed surfaces}

In this section we collect some preliminary results on surface braid groups that are needed in the sequel of the work.

\begin{definition-conf}
Let $X$ be a topological space. The $k$th \emph{ordered configuration space} of $X$ is defined as
\begin{equation} \label{eq:n-conf}
\mathrm{Conf}_k (X):=\{(x_1,\ldots, x_k) \in X^k \, | \, \ x_i \neq x_j \;\; \textrm{for all} \; \;i \neq j\}, 
\end{equation}
namely $\mathrm{Conf}_k (X) = X^k - \Delta$, where $\Delta$ is the big diagonal. 
 
The quotient of $\mathrm{Conf}_k (X)$ by the natural free action of the symmetric group $\mathfrak{S}_k$ is called the $k$th \emph{unordered configuration space} of $X$, and it is denoted by $\mathrm{UConf}_k (X)$.  

Then $\mathrm{UConf}_k (X) = \mathrm{Sym}^k(X) - \delta$, where $\delta$ denotes the image of $\Delta$ in the symmetric product. 
\end{definition-conf}

\begin{remark}
If $X$ is a smooth, compact, $n$-dimensional manifold, then both the configuration spaces $\mathrm{Conf}_k (X)$ and $\mathrm{UConf}_k (X)$ are smooth, open, $kn$-dimensional manifolds.  
\end{remark}

\begin{definition-braid}
Let $\Sigma_g$ be a closed topological surface of genus $g$, and let  $\mathscr{P} = \{p_1, \ldots, p_k\} \subset \Sigma_g$ be a set of $k$ distinct points. A \emph{geometric braid} on $\Sigma_g$ based at $\mathscr{P}$ (also called a \emph{braid on} $k$ \emph{strings}) is a $k$-ple $(\alpha_1, \ldots, \alpha_k)$ of paths $\alpha_i \colon [0, \, 1] \to \Sigma_g$ such that 
\begin{itemize}
\item $\alpha_i(0) = p_i, \quad i=1, \ldots, k$;
\item $\alpha_i(1) \in \mathscr{P}, \quad i=1, \ldots, k$;
\item the points $\alpha_1(t), \ldots, \alpha_k(t) \in \Sigma_g$ are pairwise distinct for all $t \in [0, \, 1]$.
\end{itemize}     
A geometric braid such that $\alpha_i(0)=\alpha_i(1)$ for all $i \in \{1, \ldots, k\}$ is called a \emph{pure geometric braid}.
\end{definition-braid}
\begin{figure}[h]
\centering
\begin{tikzpicture}
\braid[number of strands=3, style strands={1}{red}, style strands={2}{blue}, style strands={3}{olive}, line width=1.5pt]
(braid) a_1 a_2 a_1^{-1};
\end{tikzpicture}
\medskip
\caption{A non-pure braid on $3$ strings} \label{fig:braids}
\end{figure}
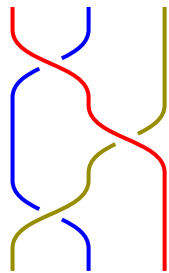
The \emph{braid group} on $k$ strings on $\Sigma_g$ is the group $\mathsf{B}_{k}(\Sigma_g)$ whose elements are the braids  based at $\mathscr{P}$ and whose operation is the usual product of paths, up to homotopies among braids. The \emph{pure braid group} is the subgroup $\mathsf{P}_{k}(\Sigma_g)$ of $\mathsf{B}_{k}(\Sigma_g)$ given by the homotopy classes of pure braids. It can be shown that $\mathsf{B}_{k}(\Sigma_g)$ and $\mathsf{P}_{k}(\Sigma_g)$ do not depend on the choice of the set $\mathscr{P}$, and that there is a short exact sequence of groups
\begin{equation} \label{eq:pure-nonpure}
1 \to \mathsf{P}_{k}(\Sigma_g)\to \mathsf{B}_{k}(\Sigma_g)\to \mathfrak{S}_k \to 1.
\end{equation}
Moreover, there are isomorphisms 
\begin{equation} \label{eq:iso-braids}
\mathsf{P}_{k}(\Sigma_g) \simeq \pi_1(\mathrm{Conf}_k (\Sigma_g)), \quad \mathsf{B}_{k}(\Sigma_g) \simeq \pi_1(\mathrm{UConf}_k (\Sigma_g)),
\end{equation}
so that we can interpret \eqref{eq:pure-nonpure} as the short exact sequence of fundamental groups induced by the $\mathfrak{S}_k$-covering 
\begin{equation}
\mathrm{Conf}_k (\Sigma_g) \to \mathrm{UConf}_k (\Sigma_g). 
\end{equation}
Braid groups are an important and flexible tool used in several areas of science, such as 
\begin{itemize}
\item[-] Knot Theory (Alexander's theorem)
\item[-] Mathematical Physics (Yang-Baxter's equation)
\item[-] Mechanical Engineering (robot motion planning)
\item[-] Algebraic Geometry (monodromy invariants).
\end{itemize}
We will focus on the last topic, explaining
how the representation monodromy of surface braid groups onto the symmetric group can be used to produce interesting examples of projective surfaces defined over the field of complex numbers. \\

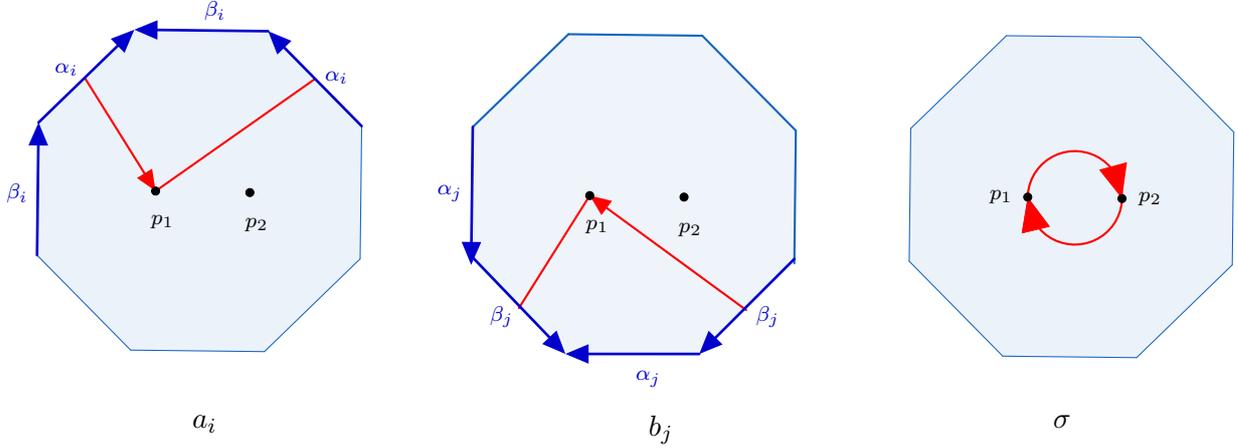
\begin{figure}[H]
 \begin{tikzpicture}[line cap=round,line join=round,>=triangle 45,x=1cm,y=1cm]
\clip(-6.17,-1.27) rectangle (15.43,5.43);
\fill[line width=0.4pt,color=rvwvcq,fill=rvwvcq,fill opacity=0.08] (-2.5,4.61) -- (-4.26,4.63) -- (-5.518650070512054,3.3996342007354072) -- (-5.538650070512054,1.6396342007354078) -- (-4.3082842712474605,0.38098413022335364) -- (-2.548284271247461,0.3609841302233536) -- (-1.289634200735407,1.5913499294879463) -- (-1.269634200735407,3.351349929487946) -- cycle;
\fill[line width=0.8pt,color=rvwvcq,fill=rvwvcq,fill opacity=0.07] (3.21,4.55) -- (1.45,4.57) -- (0.1913499294879455,3.3396342007354085) -- (0.17134992948794458,1.579634200735409) -- (1.4017157287525365,0.32098413022335403) -- (3.161715728752536,0.30098413022335313) -- (4.420365799264591,1.531349929487945) -- (4.440365799264592,3.2913499294879447) -- cycle;
\fill[line width=0.4pt,color=rvwvcq,fill=rvwvcq,fill opacity=0.08] (8.97,4.53) -- (7.21,4.55) -- (5.951349929487945,3.3196342007354067) -- (5.931349929487945,1.559634200735406) -- (7.161715728752538,0.3009841302233509) -- (8.921715728752538,0.2809841302233509) -- (10.180365799264594,1.511349929487944) -- (10.200365799264593,3.2713499294879446) -- cycle;
\draw [line width=0.4pt,color=rvwvcq] (-5.538650070512054,1.6396342007354078)-- (-4.3082842712474605,0.38098413022335364);
\draw [line width=0.4pt,color=rvwvcq] (-4.3082842712474605,0.38098413022335364)-- (-2.548284271247461,0.3609841302233536);
\draw [line width=0.4pt,color=rvwvcq] (-2.548284271247461,0.3609841302233536)-- (-1.289634200735407,1.5913499294879463);
\draw [line width=0.4pt,color=rvwvcq] (-1.289634200735407,1.5913499294879463)-- (-1.269634200735407,3.351349929487946);
\draw [line width=0.8pt,color=rvwvcq] (3.21,4.55)-- (1.45,4.57);
\draw [line width=0.8pt,color=rvwvcq] (1.45,4.57)-- (0.1913499294879455,3.3396342007354085);
\draw [line width=0.8pt,color=rvwvcq] (4.420365799264591,1.531349929487945)-- (4.440365799264592,3.2913499294879447);
\draw [line width=0.8pt,color=rvwvcq] (4.440365799264592,3.2913499294879447)-- (3.21,4.55);
\draw [line width=0.4pt,color=rvwvcq] (8.97,4.53)-- (7.21,4.55);
\draw [line width=0.4pt,color=rvwvcq] (7.21,4.55)-- (5.951349929487945,3.3196342007354067);
\draw [line width=0.4pt,color=rvwvcq] (5.951349929487945,3.3196342007354067)-- (5.931349929487945,1.559634200735406);
\draw [line width=0.4pt,color=rvwvcq] (5.931349929487945,1.559634200735406)-- (7.161715728752538,0.3009841302233509);
\draw [line width=0.4pt,color=rvwvcq] (7.161715728752538,0.3009841302233509)-- (8.921715728752538,0.2809841302233509);
\draw [line width=0.4pt,color=rvwvcq] (8.921715728752538,0.2809841302233509)-- (10.180365799264594,1.511349929487944);
\draw [line width=0.4pt,color=rvwvcq] (10.180365799264594,1.511349929487944)-- (10.200365799264593,3.2713499294879446);
\draw [line width=0.4pt,color=rvwvcq] (10.200365799264593,3.2713499294879446)-- (8.97,4.53);
\draw [line width=0.8pt,color=ffqqqq] (-3.98,2.49)-- (-1.8820135571336331,3.9778069722401543);
\draw [->,line width=0.8pt,color=ffqqqq] (-4.912072950290508,3.992580374435121) -- (-3.98,2.49);
\draw [line width=0.8pt,color=ffqqqq] (1.73,2.43)-- (0.7983789804885912,0.9381906887922181);
\draw [->,line width=0.8pt,color=ffqqqq] (3.787196012305263,0.9124086770670292) -- (1.73,2.43);
\draw [->,line width=1pt,color=qqqqcc] (-5.5386500705120545,1.6396342007354092) -- (-5.518650070512055,3.399634200735406);
\draw [->,line width=1pt,color=qqqqcc] (-5.518650070512055,3.399634200735406) -- (-4.26,4.63);
\draw [->,line width=1pt,color=qqqqcc] (-1.269634200735407,3.3513499294879456) -- (-2.5,4.61);
\draw [->,line width=1pt,color=qqqqcc] (-2.5,4.61) -- (-4.26,4.63);
\draw [->,line width=1pt,color=qqqqcc] (0.1913499294879455,3.3396342007354085) -- (0.17,1.51);
\draw [->,line width=1pt,color=qqqqcc] (0.17139736171360578,1.6297503150282006) -- (1.41,0.33);
\draw [->,line width=1pt,color=qqqqcc] (3.17,0.33) -- (1.41,0.33);
\draw [->,line width=1pt,color=qqqqcc] (4.421030852899942,1.5898746493988644) -- (3.17,0.33);
\draw [shift={(8.11,2.4)},line width=0.8pt,color=ffqqqq]  plot[domain=-0.016127633843636246:3.125465019746157,variable=\t]({1*0.6200806399171007*cos(\t r)+0*0.6200806399171007*sin(\t r)},{0*0.6200806399171007*cos(\t r)+1*0.6200806399171007*sin(\t r)});
\draw [->,line width=2pt,color=ffqqqq] (8.639103545440543,2.723341055546999) -- (8.73,2.39);
\draw [shift={(8.11,2.4)},line width=0.8pt,color=ffqqqq]  plot[domain=3.125465019746157:6.26705767333595,variable=\t]({1*0.6200806399170989*cos(\t r)+0*0.6200806399170989*sin(\t r)},{0*0.6200806399170989*cos(\t r)+1*0.6200806399170989*sin(\t r)});
\draw [->,line width=2pt,color=ffqqqq] (7.563712410516943,2.106623331573913) -- (7.49,2.41);
\draw [->,line width=2pt,color=ffqqqq] (7.563712410516943,2.106623331573913) -- (7.49,2.41);
\draw (-3.63,-0.35) node[anchor=north west] {$a_i$};
\draw (2.37,-0.35) node[anchor=north west] {$b_j$};
\draw (7.69,-0.35) node[anchor=north west] {$\sigma$};
\begin{scriptsize}
\draw [fill=black] (-3.98,2.49) circle (1.5pt); \draw[color=black]
(-3.89,2.08) node {$p_1$}; \draw [fill=black] (-2.74,2.47)
circle(1.5pt); \draw[color=black] (-2.65,2.04) node {$p_2$}; \draw
[fill=black] (1.73,2.43) circle (1.5pt); \draw[color=black]
(1.83,2.02) node {$p_1$}; \draw [fill=black] (2.97,2.41) circle
(1.5pt); \draw[color=black] (3.05,1.98) node {$p_2$}; \draw
[fill=black] (7.49,2.41) circle (1.5pt); \draw[color=black] (7.14,2.41) node {$p_1$};  \draw [fill=black]
(8.73,2.39) circle (1.5pt); \draw[color=black] (9.10,2.39) node {$p_2$}; 
\draw[color=qqqqcc] (-5.80,2.46) node
{$\beta_i$}; 
\draw[color=qqqqcc] (-5.15,4.1) node {$\alpha_i$};
\draw[color=qqqqcc] (-1.6,4) node {$\alpha_i$};
\draw[color=qqqqcc] (-3.20,4.88) node {$\beta_i$};
\draw[color=qqqqcc] (-0.11,2.46) node {$\alpha_j$};
\draw[color=qqqqcc] (0.57,0.82) node {$\beta_j$};
\draw[color=qqqqcc] (2.50,0) node {$\alpha_j$}; \draw[color=qqqqcc]
(4.07,0.82) node {$\beta_j$};
\end{scriptsize}
\end{tikzpicture}
\medskip
\caption{Generators of $\mathsf{B}_2(\Sigma_2)$}\label{fig:2}
\end{figure}

We are primarily interested in the case $g=k=2$. In that case, a simple presentation for the braid group is provided by the following result.
\begin{bellingeri} 
The braid group $\mathsf {B}_2(\Sigma_2)$ can be generated by five elements $a_1, \, a_2, \, b_1, \, b_2, \, \sigma,$
subject to the eleven relations below$:$
\begin{equation} \label{eq:relations}
\begin{split}
(R2) \quad &  \sigma^{-1} a_1 \sigma^{-1} a_1= a_1 \sigma^{-1} a_1 \sigma^{-1} \\ &  \sigma^{-1} a_2 \sigma^{-1} a_2= a_2 \sigma^{-1} a_2 \sigma^{-1} \\ &
\sigma^{-1} b_1 \sigma^{-1} b_1 = b_1 \sigma^{-1} b_1 \sigma^{-1} \\ & \sigma^{-1} b_2 \sigma^{-1} b_2 = b_2 \sigma^{-1} b_2 \sigma^{-1}\\ 
&  \\
(R3) \quad & \sigma^{-1} a_1 \sigma a_2 = a_2 \sigma^{-1} a_1 \sigma \\ &  \sigma^{-1} b_1 \sigma b_2 = b_2 \sigma^{-1} b_1 \sigma \\
& \sigma^{-1} a_1 \sigma b_2 = b_2 \sigma^{-1} a_1 \sigma \\
& \sigma^{-1} b_1 \sigma a_2 = a_2 \sigma^{-1} b_1 \sigma \\
 & \\
(R4) \quad & \sigma^{-1} a_1 \sigma^{-1} b_1 = b_1 \sigma^{-1} a_1 \sigma \\
 & \sigma^{-1} a_2 \sigma^{-1} b_2 = b_2 \sigma^{-1} a_2 \sigma \\
 & \\
 (TR) \quad &  [a_1, \, b_1^{-1}] [a_2, \, b_2^{-1}]= \sigma^2.  
\end{split}
\end{equation}
\end{bellingeri}
Here the $a_i$ and the $b_i$ are pure braids coming from the representation of the topological surface $\Sigma_2$ as a polygon of $8$ sides with the standard identification of the edges, whereas $\sigma$ is a non-pure braid exchanging the two points $p_1$, $p_2 \in \Sigma_2$. These braids are depicted in Figure \ref{fig:2}; note that, both in the cases of $a_i$ and $b_j$, the only non-trivial string is the first one, which goes through the wall $\alpha_i$, respectively the wall $\beta_j$. 
Regarding the generator $\sigma$, in terms of the isomorphism 
\begin{equation}
\mathsf{B}_{2}(\Sigma_2) \simeq \pi_1(\mathrm{UConf}_2 (\Sigma_2)) = \pi_1(\mathrm{Sym}^2(\Sigma_2)-\delta), 
\end{equation}
it corresponds to the homotopy class in $\mathrm{UConf}_2(\Sigma_2)$ of a topological loop in $\mathrm{Sym^2}(\Sigma_2)$ that ``winds once around the diagonal $\delta$".

\section{Finite coverings and monodromy representations}

Let us recall now the classification of branched coverings 
$f \colon X \to Y$ of complex, projective varieties via the classification of monodromy representations of the fundamental group $\pi_1(Y-B)$, where $B \subset Y$ is the branch locus of $f$. The main technical tools needed are the \emph{Grauert-Remmert extension theorem} and the GAGA \emph{principle}, that we recall below. 

\begin{GR}
Let $Y$ be a normal analytic space over $\mathbb{C}$ and $Z \subset Y$ a closed analytic subspace such that $U=Y - Z$ is dense in $Y$. Then any finite, analytic, unramified covering $f^{\circ} \colon V \to U$ can be extended to a normal, analytic, finite covering $f \colon X \to Y$, branched at most over $Z$. Furthermore, such an extension is unique up to analytic isomorphisms.   
\end{GR}

\begin{GAGA}
Let $X$, $Y$ be projective varieties over $\mathbb{C}$, and $X^{\rm an}$, $Y^{\rm an}$ the underlying complex analytic spaces. Then
\begin{itemize}
\item[$\boldsymbol{(i)}$] every analytic map $X^{\rm an} \to Y^{\rm an}$ is algebraic$;$
\item[$\boldsymbol{(ii)}$] every coherent analytic sheaf on $X^{\rm an}$ is algebraic, and its algebraic cohomology coincides with its analytic one.
\end{itemize} 
\end{GAGA}
From this, we deduce the following important consequences.
\begin{proj-ext}  
Let $Y$ be a smooth, projective variety over $\mathbb{C}$ and $Z \subset Y$ be a smooth, irreducible divisor. Set $U=Y - Z$. Then any finite, unramified analytic covering $f^{\circ} \colon V \to U$ can be extended in a unique way to a finite covering $f \colon X \to Y,$ branched at most over $Z$. \vskip 3mm

Moreover, there exists on $X$ a unique structure of smooth projective variety that makes $f$ an algebraic finite covering. 
\end{proj-ext}

\begin{monodromy} 
Let $Y$ be a smooth projective variety over $\mathbb{C}$ and $Z \subset Y$ be a smooth, irreducible divisor. Then isomorphism classes of connected coverings of degree $n$
\begin{equation*}
f \colon X \to Y,
\end{equation*}
branched at most over $Z$, are in bijection to group homomorphisms with transitive image
\begin{equation} 
\varphi \colon \pi_1(Y - Z) \to \mathfrak{S}_n,
\end{equation}
up to conjugacy in $\mathfrak{S}_n$. Furthermore, $f$ is a Galois covering if and only if the subgroup $\mathrm{im}\, \varphi$ of $\mathfrak{S}_n$ has order $n$, and in this case $\mathrm{im}\, \varphi$ is isomorphic to the Galois group of $f$.
\end{monodromy}
The group homomorphism 
\begin{equation} 
\varphi \colon \pi_1(Y - Z) \to \mathfrak{S}_n,
\end{equation}
is called the \emph{monodromy representation} of the covering $f \colon X \to Y$, and its image $\mathrm{im}\, \varphi \subseteq \mathfrak{S}_n$ is called the \emph{monodromy group} of $f$. The last corollary implies that, if $f$ is a Galois covering, then the monodromy group of $f$ is isomorphic to its Galois group (coinciding with the group $D(X/Y)$ of deck transformations of the covering).

\section{Generic coverings of $\mathrm{UConf}_2(\Sigma_2)$}

We now apply the previous theory in the special case
\begin{equation}
Y = \mathrm{UConf(\Sigma_2)}=\mathrm{Sym}^2(\Sigma_2)-\delta, \quad  Z  = \delta.
\end{equation}
We consider $\Sigma_2$ as a \emph{compact Riemann surface}, namely we fix a complex structure on it. Then the Abel-Jacobi map
\begin{equation}
\pi \colon \mathrm{Sym}^2(\Sigma_2) \to J(\Sigma_2)
\end{equation}
is a birational morphism onto the the Abelian surface $J(\Sigma_2)$, more precisely it is the blow-down of the unique rational curve $E \subset \mathrm{Sym}^2(\Sigma_2)$, namely the $(-1)$-curve given by the graph of the hyperelliptic involution on $\Sigma_2$. 

We have $\delta E=6$, because the curve $E$ intersects the diagonal $\delta$ transversally at the six points corresponding to the six Weierstrass points of $\Sigma_2$. Writing $\Theta$ for the numerical class of a Theta divisor in $J(\Sigma_2)$, it follows that the image $D:=\pi_*\delta \subset J(\Sigma_2)$ is an irreducible curve with an ordinary sextuple point and no other singularities, whose numerical class is $4 \Theta$.  

\begin{definition-gen-cov}
Let $f \colon S \to \mathrm{Sym}^2(\Sigma_2)$ be a connected covering of degree $n$
branched over the diagonal $\delta$, with ramification divisor $R \subset S$. Then $f$ is called \emph{generic} if  
\begin{equation*}
f^* \delta= 2R + R_0,
\end{equation*}
where the restriction $\left.f\right|_{R} \colon R \to \delta$ is an isomorphism and $R_0$ is an effective divisor over which $f$ is not ramified. 
\end{definition-gen-cov}
Generic coverings are never Galois, unless $n=2$, in which case $R_0$ is empty and  $f^* \delta = 2R$. 
Since $\delta$ is smooth, the genericity condition in the previous definition  is equivalent to requiring that the fibre of $f$ over every point of $\delta$ has cardinality $n-1$.  

\begin{figure}[H]
\centering
\begin{tikzpicture}[xscale=-1,yscale=-0.25,inner sep=0.7mm,place/.style={circle,draw=black!100,fill=black!100,thick}] 
\draw (-0.7,-0.5) rectangle (0.6,5);

\draw[red,rotate=92,x=6.28ex,y=1ex] (0.9,-0.85) cos (1,0) sin (1.25,1) cos (1.5,0) sin (1.75,-1) cos (2,0) sin (2.25,1) cos (2.5,0) sin (2.75,-1) cos (3,0) sin (3.25,1) cos (3.5,0) sin (3.6,-0.85);
\draw[rotate=92] (0.5,-0.025) .. controls (1.75,0.035) and (2.75,0.035) .. (4,-0.025);

\draw (-6.7,-0.5) rectangle (-5.4,5);

\draw[red,rotate=92,x=6.28ex,y=1ex,xshift=6,yshift=170] (0.9,-0.85) cos (1,0) sin (1.25,1) cos (1.5,0) sin (1.75,-1) cos (2,0) sin (2.25,1) cos (2.5,0) sin (2.75,-1) cos (3,0) sin (3.25,1) cos (3.5,0) sin (3.6,-0.85);
\draw[rotate=92,xshift=6,yshift=170] (0.5,-0.025) .. controls (1.75,0.035) and (2.75,0.035) .. (4,-0.025);

\draw (-6.7,15.5) rectangle (-5.4,21);

\draw[red,xscale=-1,yscale=-4] (6.05,-4.55) node(A0) [place,scale=0.2]{} to [in=5,out=55,looseness=8mm,loop] () to [in=65,out=115,looseness=8mm,loop] () to [in=125,out=175,looseness=8mm,loop] () to [in=185,out=235,looseness=8mm,loop] () to [in=245,out=295,looseness=8mm,loop] () to [in=305,out=355,looseness=8mm,loop] ();

\draw[xscale=-1,yscale=-4] (0.2,-0.1) node(B0) []{\textcolor{red}{\footnotesize{$R$}}};
\draw[xscale=-1,yscale=-4] (0.25,-1.05) node(0B) []{\footnotesize{$Z$}};

\draw[xscale=-1,yscale=-4] (6.25,-0.1) node(B1) []{\textcolor{red}{\footnotesize{$\delta$}}};
\draw[xscale=-1,yscale=-4] (6.25,-1.05) node(1B) []{\footnotesize{$E$}};

\draw[xscale=-1,yscale=-4] (6.4,-5) node(B2) []{\textcolor{red}{\footnotesize{$D$}}};

\draw[xscale=-1,yscale=-4] (-0.9,-0.125) node(C0) []{$S$};
\draw[xscale=-1,yscale=-4] (7,-0.125) node(C1) []{$\quad \quad \quad \;  \mathrm{Sym}^2(\Sigma_2)$};
\draw[xscale=-1,yscale=-4] (7,-5.00) node(C2) []{$\quad \quad J(\Sigma_2)$};

\draw[->] (-0.9,2.25) -- (-5.2,2.25) node[midway,above] {$f$};

\draw[->] (-6.05,5.8) -- (-6.05,14.7) node[midway,right] {$\pi$};

\draw[->] (-0.9,4.8) -- (-5.2,15.5) node[midway,below=3pt] {$\alpha$}; 
\end{tikzpicture}
\medskip
\caption{A generic covering of $\mathrm{Sym}^2(\Sigma_2)$ branched over $\delta$ }\label{fig:3}
\end{figure}
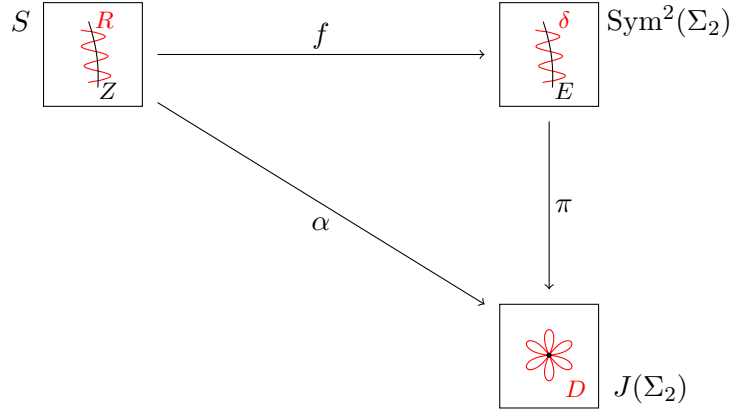

Setting $\alpha:=\pi \circ f$, the case where the curve  $Z=f^*(E)$ is irreducible is illustrated in Figure \ref{fig:3}.

From now on, by \emph{surface} we mean a smooth, complex, projective variety $S$ with $\dim_{\mathbb{C}}(S)=2$. For such a surface
 \begin{itemize}
\item $K_S=\wedge^2 \Omega_S^1$ denotes the \emph{canonical line bundle} 
\item $p_g(S)=h^0(S, \, K_S)$ is the \emph{geometric genus}
\item $q(S)=h^1(S, \, K_S)$ is the \emph{irregularity}
\item $\chi(\mathcal{O}_S)=1-q(S)+p_g(S)$ is the \emph{holomorphic Euler-Poincar\'e characteristic}.
\end{itemize}
The zero locus of a meromorphic section of $K_S$ defines a class in $H_2(S, \, \mathbb{Z})$, whose Poincar\'e 
dual $[K_S] \in H^2(S, \, \mathbb{Z})$ is called the \emph{canonical class} of $S$. Its self-intersection is the integer number defined by the cup product
\begin{equation}
K_S^2=[K_S] \cup [K_S] \in H^4(S, \, \mathbb{Z}) \simeq \mathbb{Z}.
\end{equation}
Given a surface $S$, we can define some important meromorphic maps on it: the \emph{pluricanonical maps} and the \emph{Albanese map}. 
\begin{definition-pluri}
Set $N(r):=\dim H^0(S, \, K_S^{\otimes r})-1$ and let $\{ \sigma_0, \ldots, \sigma_{N(r)} \}$ be a basis for $H^0(S, \, K_S^{\otimes r})$. Then the $r$-\emph{th} \emph{pluricanonical map} of $S$ is  the rational map
\begin{equation}
\psi_r \colon S \to \mathbb{P}^{N(r)}(\mathbb{C}), \quad x \mapsto [\sigma_0(x): \ldots : \sigma_{N(r)}(x)].
\end{equation}
We say that $S$ is \emph{of general type} if the image of $\psi_r$ is a surface for $r$ large enough (i.e., if $\psi_r$ is generically finite onto its image for $r$ large enough).
\end{definition-pluri}
\begin{definition-alb}
The \emph{Albanese map} of $S$ is the rational morphism
\begin{equation}
a_S \colon S \to \mathrm{Alb}(S):=H^0(S, \, \Omega^1_S)^*/H_1(S, \, \mathbb{Z}), 
\end{equation}
defined by the integration of global, holomorphic $1$-forms on $S$ (this is a generalization of the Abel-Jacobi map $C \to J(C)$, sending a smooth complex curve into its Jacobian). Note that $\mathrm{Alb}(S)$ is a complex torus (actually, an Abelian variety) of dimension $q(S)$.

We say that $S$ is \emph{of maximal Albanese dimension} if the image of its Albanese map is a surface 
(i.e., if $a_S$ is generically finite onto its image); note that this condition implies $q(S) \geq 2$.
\end{definition-alb}
Surfaces of general type satisfy $\chi(\mathcal{O}_S) \geq 1$, and those with $\chi(\mathcal{O}_S)=1$ are usually difficult to construct. This explains the relevance of the following result.
\begin{theoremP1}
Let $f \colon S \to \mathrm{Sym}^2(\Sigma_2)$ be a generic covering of degree $n$ and whose branch locus is the diagonal $\delta$. Then $S$ is a surface of maximal Albanese dimension with
\begin{equation*}
\chi(\mathcal{O}_S)=1, \quad K_S^2= 10-n. 
\end{equation*}
Moreover, if $ 2 \leq n \leq 9$ then $S$ is of general type.
\end{theoremP1}
Our aim is now to use the theory developed before in order to construct generic coverings $f \colon S \to \mathrm{Sym}^2(\Sigma_2)$.

\begin{definition-gen-mod} 
A \emph{generic monodromy representation} of the braid group $\mathsf{B}_2(\Sigma_2)$ is a group homomorphism 
\begin{equation}
\varphi \colon \mathsf{B}(\Sigma_2) \to \mathfrak{S}_n
\end{equation}
with transitive image and such that $\varphi(\sigma)$ is a \emph{transposition}.
\end{definition-gen-mod}
By the previous discussion, since $\mathsf{B}_2(\Sigma_2) \simeq \pi_1(\mathrm{Sym}^2(\Sigma_2) - \delta)$, generic coverings and generic monodromy representations are related by the following
\begin{theoremP2} 
Isomorphism classes of generic coverings of degree $n$ 
\begin{equation}
f \colon S \to \mathrm{Sym}^2(\Sigma_2),
\end{equation}
with branched locus $\delta$, are in bijective correspondence to generic monodromy representations
\begin{equation}
\varphi \colon \mathsf{B}_2(\Sigma_2) \to \mathfrak{S}_n,
\end{equation}    
up to conjugacy in $\mathfrak{S}_n$. 
For $2 \leq n \leq 9$, the number of such representations is given in the table below$:$  
\begin{table}[H]
\begin{center}
\begin{tabular}{c|c|c|c|c|c|c|c|c}
$ n$ & $2$ & $3$ & $4$ & $5$ & $6$ & $7$ & $8$ & $9$ \\
 \hline
$\mathrm{Number \; of}$ $\varphi$ & $16$ & $3 \cdot 80$ & $6 \cdot 480$ & $0$ & $15 \cdot 2880$  & $0$ & $28 \cdot 172800$ & $0$\\
\end{tabular}
\end{center}
\end{table} In particular, for $n \in \{5, \, 7, \, 9\}$ there exist no generic coverings.
\end{theoremP2}

\bigskip \bigskip

As an application of the general theory, let us finish this note by discussing in detail the cases $n=2$ and $n=3$.
\begin{example-2}
In this case we are looking for generic monodromy representations
\begin{equation}
\varphi \colon \mathsf {B}_2(\Sigma_2) \to \mathfrak{S}_2 = \{(1), \, (1 \, 2)\}.
\end{equation} 
Since $\mathsf{B}_2(\Sigma_2)$ is generated by five elements $a_1$, $a_2$, $b_1$, $b_2$, $\sigma$, and necessarily $\varphi(\sigma)=(1 \, 2)$, we see that there are $2^4=16$ possibilities for $\varphi$. 

The group $\mathfrak{S}_2$ is abelian, so there is no conjugacy relation to consider and we get sixteen isomorphism classes of double coverings $f \colon S \to \mathrm{Sym}^2(\Sigma_2)$, branched over $\delta$ and with 
\begin{equation}
\chi(\mathcal{O}_S)=1, \quad K_S^2=8.
\end{equation}
These coverings correspond to the sixteen square roots of $\delta$ in  the Picard group of $\mathrm{Sym}^2(\Sigma_2)$. One covering coincides with the natural projection $f \colon C_2 \times C_2 \to \mathrm{Sym}^2(\Sigma_2)$, in fact 
\begin{equation*}
p_g(C_2 \times C_2) = q(C_2 \times C_2) =4, \quad K_{C_2 \times C_2}=8.
\end{equation*}
The remaining fifteen coverings are surfaces of general type  with
\begin{equation}
p_g(S)=q(S)=2, \quad K_S^2=8.
\end{equation}
\end{example-2}
\begin{example-3}
In this case we are looking for generic monodromy representations
\begin{equation*}
\varphi \colon \mathsf{B}_2(\Sigma_2) \to \mathfrak{S}_3, 
\end{equation*} 
up to conjugacy in $\mathfrak{S}_3$. By hands, or by using a Computer Algebra software like $\mathsf{GAP4}$ (see \cite{5}), we find that the total  number of monodromy representations is $240$. For every such a representation we have $\mathrm{im}\, \varphi = \mathfrak{S}_3$. Moreover, each orbit for the conjugacy action of $\mathfrak{S}_3$ on the set of monodromy representations consists of six elements, and consequently the orbit set has cardinality $240/6 = 40$. 
By our last theorem, this implies that there are $40$ isomorphism classes of generic coverings $f \colon S \to \mathrm{Sym}^2(\Sigma_2)$ of degree $3$ and branched over $\delta$. For all of them, $S$ is a surface of general type with
\begin{equation*}
  p_g(S)=q(S)=2, \quad K_S^2=7 
\end{equation*}
and its Albanese map $\alpha \colon S \to J(\Sigma_2)$ is a generically finite covering of degree $3$.

These surfaces were previously studied by R. Pignatelli and the Author in \cite{6}, by using completely different methods. In fact, they showed that they all lie in the same deformation class, and that their moduli space is a connected, quasi-finite cover of degree $40$ of $\mathcal{M}_2$, the coarse moduli space of curves of genus $2$.
\end{example-3}

\section*{Acknowledgments}
The Author wishes to thank the organizers of Group 32 - The 32nd \emph{International Colloquium on Group Theoretical Methods in Physics}, held on Czech Technical University (Prague) on July 9-13, 2018, for the invitation and the hospitality.

\end{document}